\documentclass[12pt]{amsart}
\usepackage{amsmath,amssymb,amsthm,amscd,amsxtra,eucal,latexsym,mathrsfs}

\newcommand{\nc}{\newcommand}

\newcommand{\cA}{{\mathcal A}}
\newcommand{\cB}{{\mathcal B}}

\newcommand{\cE}{{\mathcal E}}

\newcommand{\cI}{{\mathcal I}}

\newcommand{\cO}{{\mathcal O}}
\newcommand{\cL}{{\mathcal L}}

\newcommand{\cF}{{\mathcal F}}

\newcommand{\cP}{{\mathcal P}}

\newcommand{\cV}{{\mathcal V}}
\newcommand{\cW}{{\mathcal W}}

\renewcommand{\AA}{{\mathbb A}}

\newcommand{\NN}{{\mathbb N}}
\newcommand{\ZZ}{{\mathbb Z}}
\newcommand{\QQ}{{\mathbb Q}}

\newcommand{\PP}{{\mathbb P}}

\newcommand{\EE}{{\mathbb E}}

\renewcommand{\gg}{\mathfrak{g}}  

\newcommand{\gt}{\mathfrak{t}}
\newcommand{\gr}{\mathfrak{r}}
\newcommand{\gs}{\mathfrak{s}}

\newcommand{\on}{\operatorname}

\newcommand{\Fl}{{\mathcal F}l}

\newcommand{\Rep}{{\on{Rep}}}
\newcommand{\Sch}{{\on{Sch}}}

\newcommand{\Qlb}{\mathbb{\bar Q}_\ell}
\newcommand{\Gm}{\mathbb{G}_m}

\newcommand{\Ql}{\mathbb{Q}_\ell}
\newcommand{\toup}[1]{\stackrel{#1}{\to}}

\newcommand{\hook}[1]{\stackrel{#1}{\hookrightarrow}}

\newcommand{\getsup}[1]{\stackrel{#1}{\gets}}

\newcommand{\Hom}{\on{Hom}}

\newcommand{\Ker}{\on{Ker}}

\newcommand{\Aut}{\on{Aut}}

\newcommand{\RG}{\on{R\Gamma}}

\newcommand{\Bun}{\on{Bun}}

\newcommand{\Spec}{\on{Spec}}

\newcommand{\HOM}{{{\mathcal H}om}}

\newcommand{\Gr}{\on{Gr}}

\newcommand{\pr}{\on{pr}}
\newcommand{\id}{\on{id}}

\newcommand{\QED}{$\square$} 
  
\newcommand{\iso}{{\widetilde\to}}

\newcommand{\comp}{\circ}

\newcommand{\DD}{\mathbb{D}}  
\newcommand{\D}{\on{D}}       
\newcommand{\wt}{\widetilde}
\newcommand{\ov}[1]{\overline{#1}}
\newcommand{\select}[1]{{\it{#1}}}

\renewcommand{\P}{{\on{P}}}

\newcommand{\<}{\langle}
\renewcommand{\>}{\rangle}

\newcommand{\Conv}{\on{Conv}}

\newcommand{\Lie}{\on{Lie}}

\newcommand{\Vect}{\on{Vect}}

\nc{\Perv}{\on{Perv}}

\newcommand{\cPic}{\on{{\cP}ic}}

\nc{\Gra}{\on{Gra}}
\nc{\PPerv}{\on{{\PP}erv}}

\nc{\oX}{\overset{\scriptscriptstyle\circ}{X}}
\nc{\gRes}{\on{gRes}}
\nc{\ocL}{\overset{\circ}{\cL}}
\nc{\RHom}{\on{RHom}}
\nc{\slLie}{\mathfrak{sl}}

\newtheorem{Lm}{Lemma}[section]
\newtheorem{Th}{Theorem}[section]

\theoremstyle{remark}
\newtheorem{Rem}{Remark}[section]

\theoremstyle{definition}

\newenvironment{Prf}{\par\noindent {\it Proof }}{\QED}

\begin{document}

\title{Twisted geometric Satake equivalence: reductive case}
\author{Sergey Lysenko}
\begin{abstract} 
In this paper we extend the twisted Satake equivalence established in \cite{FL} for almost simple groups to the case of split reductive groups.
\end{abstract} 
\maketitle

\section{Introduction} Let $G$ be a connected reductive group over an algebraically closed field. Brylinski-Deligne have developed the theory of central extensions of $G$ by $K_2$. According to Weissman \cite{W}, this is a natural framework for the representation theory of metaplectic groups over local and global fields (allowing to formulate a conjectural extension of the Langlands program for metaplectic groups). The geometric Langlands program naturally extends to this setting. As a step in this direction, in this paper we extend the twisted Satake equivalence established in \cite{FL} for almost simple groups to the case of reductive groups. Our input data model an extension of $G$ by $K_2$ (and cover all the isomorphism classes of such extensions). 

 In \cite{R} a similar problem was considered under an additional assumption that $k$ is of characteristic zero. However, \cite{R} is written rather on a physical than mathematical level, many definitions are not formulated precisely (even the main definition of a symmetric factorizable gerb), it also contains mistakes. In particular, (\cite{R}, Proposition~II.3.6) is wrong as stated. Our purpose is to provide a reliable reference in particular for our further works \cite{L1, L2, L3}. Our main result is Theorem~\ref{Th_one} in Section~\ref{Section_MR}. 

\section{Main result}

\subsection{Notations} 
\label{Section_notations}
Let $k$ be an algebraically closed field. Let $G$ be a split reductive group over $k$, $T\subset B\subset G$ be a maximal torus and a Borel subgroup. Let $\Lambda$ (resp., $\check{\Lambda}$) denote the coweights (resp., weights) lattice of $T$. Let $W$ denote the Weyl group of $(T,G)$. Set $\cO=k[[t]]\subset F=k((t))$. As in (\cite{L}, Section~3.2), we denote by $\cE^s(T)$ the category of pairs: a symmetric bilinear form $\kappa: \Lambda\otimes\Lambda\to\ZZ$, and a central super extension $1\to k^*\to \tilde\Lambda^s\to\Lambda\to 1$ whose commutator is $(\gamma_1, \gamma_2)_c=(-1)^{\kappa(\gamma_1,\gamma_2)}$. 

 Let $X$ be a smooth projective connected curve over $k$. Write $\Omega$ for the canonical line bundle on $X$. Fix once and for all a square root $\Omega^{\frac{1}{2}}$ of $\Omega$. 

Let $\cP^{\theta}(X,\Lambda)$ denote the category of theta-data (\cite{BD2}, Section~3.10.3). Recall the functor $\cE^s(T)\to \cP^{\theta}(X,\Lambda)$ defined in (\cite{L}, Lemma~4.1). Let $(\kappa, \tilde\Lambda^s)\in \cE^s(T)$, so for $\gamma\in\Lambda$ we are given a super line $\epsilon^{\gamma}$ and isomorphisms $c^{\gamma_1,\gamma_2}: \epsilon^{\gamma_1}\otimes \epsilon^{\gamma_2}\,\iso\, \epsilon^{\gamma_1+\gamma_2}$. For $\gamma\in\Lambda$ let $\lambda^{\gamma}=(\Omega^{\frac{1}{2}})^{\otimes -\kappa(\gamma,\gamma)}\otimes\epsilon^{\gamma}$. For the evident isomorphisms $'c^{\gamma_1,\gamma_2}: \lambda^{\gamma_1}\otimes \lambda^{\gamma_2}\,\iso\, \lambda^{\gamma_1+\gamma_2}\otimes\Omega^{\kappa(\gamma_1,\gamma_2)}$ then $(\kappa, \lambda, {'c})\in \cP^{\theta}(X,\Lambda)$. This is the image of $(\kappa, \tilde\Lambda^s)$ by the above functor.

 Let $\Sch/k$ denote the category of $k$-schemes of finite type with Zarisky topology. The $n$-th Quillen K-theory group of a scheme form a presheaf on $\Sch/k$ as the scheme varies. As in \cite{BD}, $K_n$ denotes the associated sheaf on $\Sch/k$ for the Zariski topology. 

 Denote by $\Vect$ the tensor category of vector spaces. Pick a prime $\ell$ invertible in $k$, write $\Qlb$ for the algebraic closure of $\Ql$. We work with (perverse) $\Qlb$-sheaves for \'etale topology.
 
\subsection{Motivation}
According to Weissman \cite{W}, the metaplectic input datum is an integer $n\ge 1$ and an extension $1\to K_2\to E\to G\to 1$ over $k$ as in \cite{BD}. It gives rise
to a $W$-invariant quadratic form $Q: \Lambda\to \ZZ$, for which we get the corresponding even symmetric bilinear form $\kappa: \Lambda\otimes\Lambda\to\ZZ$ given by $\kappa(x_1,x_2)=Q(x_1+x_2)-Q(x_1)-Q(x_2)$, $x_i\in\Lambda$.

 The extension $E$ yields an extension
$$
1\to K_2(F)\to E(F)\to G(F)\to 1
$$
The tame symbol gives a map $(\cdot, \cdot)_{st}: K_2(F)\to k^*$. The push-out by this map is an extension 
$$
1\to k^*\to \EE(k)\to G(F)\to 1
$$
It is the set of $k$-points of an extension of group ind-schemes over $k$
\begin{equation}
\label{ext_cE}
1\to\Gm\to \EE\to G(F)\to 1
\end{equation}
Assume $n\ge 1$ invertible in $k$. For a character $\zeta: \mu_n(k)\to\Qlb^*$ denote by $\cL_{\zeta}$ the corresponding Kummer sheaf on $\Gm$.

Pick an injective character $\zeta: \mu_n(k)\to\Qlb^*$. For a suitable section of (\ref{ext_cE}) over $G(\cO)$, we are interested in the category $\Perv_{G,\zeta}$ of $G(\cO)$-equivariant $\Qlb$-perverse sheaves on $\EE/G(\cO)$ with $\Gm$-monodromy $\zeta$, that is, equipped with $(\Gm, \cL_{\zeta})$-equivariant structure. One wants to equip it with a structure of a symmetric monoidal category (and actually a structure of a chiral category as in \cite{Ga2}), and prove a version of the Satake equivalence for it.

\subsubsection{} One has the exact sequence $1\to T_1\to T\to G/[G,G]\to 1$, where $T_1\subset [G,G]$ is a maximal torus. Write $\Lambda_{ab}$ (resp., $\check{\Lambda}_{ab}$) for the coweights (resp., weights) lattice of $G_{ab}=G/[G,G]$. The kernel of $\Lambda\to\Lambda_{ab}$ is the rational closure in $\Lambda$ of the coroots lattice. Let $J$ denote the set of connected components of the Dynkin diagram, $\cI_j$ denote the set of vertices of the $j$-th connected component of the Dynkin diagram, $\cI=\cup_{j\in J} \; \cI_j$ the set of vertices of the Dynkin diagram. For $i\in \cI$ let $\alpha_i$ (resp., $\check{\alpha}_i$) be the corresponding simple coroot (resp., root). One has $G_{ad}=\prod_{j\in J} G_j$, where $G_j$ is a simple group. Let $\gg_j=\Lie G_j$. 

 Write $\Lambda_{ad}$ for the coweights lattice of $G_{ad}$. Write $R_j$ (resp., $\check{R}_j$) for the set coroots (resp., roots) of $G_j$. Let $R$ (resp. $\check{R}$) denote the set of coroots (resp., roots) of $G$. For $j\in J$ let $\kappa_j: \Lambda_{ad}\otimes\Lambda_{ad}\to \ZZ$ denote the Killing form for $G_j$, that is, 
$$
\kappa_j=\sum_{\check{\alpha}\in \check{R}_j} \check{\alpha}\otimes\check{\alpha}
$$ 
Note that $\frac{\kappa_j}{2}: \Lambda_{ad}\otimes\Lambda_{ad}\to \ZZ$. We also view $\kappa_j$ if necessary as a bilinear form on $\Lambda$. 

 There is $m\in \NN$ such that $m\kappa$ is of the form
$$
\bar\kappa=-\beta-\sum_{j\in J} c_j\kappa_j
$$ 
for some $c_j\in\ZZ$ and some even symmetric bilinear form $\beta: \Lambda_{ab}\otimes\Lambda_{ab}\to\ZZ$. So, relaxing our assumption on the characteristic, the following setting is sufficient.

\subsection{Input data} 
\label{Section_Input_data}
For each $j\in J$ let $\cL_j$ be the ($\ZZ/2\ZZ$-graded purely of parity zero) line bundle on $\Gr_G$ whose fibre at $gG(\cO)$ is $\det(\gg_j(\cO): \gg_j(\cO)^g)$. Write $E^a_j$ for the punctured total space of the line bundle $\cL_j$ over $G(F)$. This is a central extension
\begin{equation}
\label{ext_E^a_j}
1\to \Gm\to E^a_j\to G(F)\to 1,
\end{equation}
here $a$ stands for `adjoint'. It splits canonically over $G(\cO)$. The commutator of (\ref{ext_E^a_j}) on $T(F)$ is given by
$$
(\lambda_1\otimes f_1, \lambda_2\otimes f_2)_c=(f_1, f_2)_{st}^{-\kappa_j(\lambda_1,\lambda_2)}
$$
for $\lambda_i\in \Lambda$, $f_i\in F^*$. Recall that for $f,g\in F^*$ the tame symbol is given by
$$
(f, g)_{st}=(-1)^{v(f)v(g)}(g^{v(f)}f^{-v(g)})(0) 
$$
Assume also given a central extension
\begin{equation}
\label{ext_E_beta_of_G_ab}
1\to \Gm\to E_{\beta}\to G_{ab}(F)\to 1
\end{equation}
in the category of group ind-schemes whose commutator $(\cdot, \cdot)_c: G_{ab}(F)\times G_{ab}(F)\to \Gm$ satisfies
$$
(\lambda_1\otimes f_1, \lambda_2\otimes f_2)_c=(f_1, f_2)_{st}^{-\beta(\lambda_1,\lambda_2)}
$$
for $\lambda_i\in \Lambda_{ab}$, $f_i\in F^*$. Here $\beta: \Lambda_{ab}\otimes\Lambda_{ab}\to\ZZ$ is an even symmetric bilinear form. This is a Heisenberg $\beta$-extension (\cite{BD2}, Definition~10.3.13). Its pull-back under $G(F)\to G_{ab}(F)$ is also denoted $E_{\beta}$ by abuse of notations. Assume also given a splitting of $E_{\beta}$ over $G_{ab}(\cO)$. 

 Let $N\ge 1$, assume $N$ invertible in $k$. Let $\zeta: \mu_N(k)\to \Qlb^*$ be an injective character. Assume given $c_j\in\ZZ$ for $j\in J$.
 
 The sum of the extensions $(E_j^a)^{c_j}$, $j\in J$ and the extension $E_{\beta}$ is an extension
\begin{equation}
\label{ext_cE_of_G(F)_by_T_J}
1\to \Gm\to \EE\to G(F)\to 1
\end{equation}
equipped with the induced section over $G(\cO)$. 
Set $\Gra_G=\EE/G(\cO)$. Let $\Perv_{G,\zeta}$ denote the category of $G(\cO)$-equivariant perverse sheaves on $\Gra_G$ with $\Gm$-monodromy $\zeta$. This means, by definition, a $(\Gm, \cL_{\zeta})$-equivariant structure. Set
$$
\PPerv_{G,\zeta}=\Perv_{G,\zeta}[-1]\subset \D(\Gra_G)
$$

 Let $\Gm$ act on $\EE$ via the homomorphism $\Gm\to\Gm$, $z\mapsto z^N$. Let $\wt\Gr_G$ denote the stack quotient of $\Gra_G$ by this action of $\Gm$. We view $\PPerv_{G,\zeta}$ as a full subcategory of the category of perverse sheaves on $\wt\Gr_G$ via the functor $K\mapsto \pr^*K$. Here $\pr: \Gra_G\to \wt\Gr_G$ is the quotient map. As in \cite{FL}, the above cohomological shift is a way to avoid some sign problems. 
 
 Let us make a stronger assumption that we are given a central extension
\begin{equation}
\label{ext_G_ab_by_K_2}
1\to K_2\to \cV_{\beta}\to G_{ab}\to 1
\end{equation} 
as in \cite{BD} such that passing to $F$-points and further taking the push-out by the tame symbol $K_2(F)\to \Gm$ yields the extension  (\ref{ext_E_beta_of_G_ab}). Recall  that on the level of ind-schemes the tame symbol map 
\begin{equation}
\label{tame_symbol}
(\cdot, \cdot)_{st}: F^*\times F^*\to\Gm
\end{equation}
is defined in \cite{CC}, see also (\cite{BBE}, Sections~3.1-3.3). Assume that the splitting of (\ref{ext_E_beta_of_G_ab}) over $G(\cO)$ is the following one. The composition $K_2(\cO)\to K_2(F)$ with the tame symbol map factors through $1\hook{} \Gm$, hence a canonical section $G_{ab}(\cO)\to E_{\beta}$ of (\ref{ext_E_beta_of_G_ab}). Denote by 
\begin{equation}
\label{ext_Lambda_ab_by_Gm}
1\to \Gm\to V_{\beta}\to \Lambda_{ab}\to 1
\end{equation}
the pull-back of (\ref{ext_E_beta_of_G_ab}) by $\Lambda_{ab}\to G_{ab}(F)$, $\lambda\mapsto t^{\lambda}$. This is the central extension over $k$ corresponding to (\ref{ext_G_ab_by_K_2}) by the Brylinski-Deligne classification \cite{BD}. The extension (\ref{ext_Lambda_ab_by_Gm}) is given by a line $\epsilon^{\gamma}$ (of parity zero as $\ZZ/2\ZZ$-graded) for each $\gamma\in\Lambda_{ab}$ together with isomorphisms 
$$
c^{\gamma_1,\gamma_2}: \epsilon^{\gamma_1}\otimes \epsilon^{\gamma_2}\,\iso\, \epsilon^{\gamma_1+\gamma_2}
$$ 
for $\gamma_i\in\Lambda_{ab}$ subject to the conditions in the definition of $\cE^s(T)$ (\cite{L}, Section~3.2.1). Let 
\begin{equation}
\label{ext_V_EE_of_Lambda}
1\to\Gm\to V_{\EE}\to \Lambda\to 1
\end{equation}
be the pull-back of (\ref{ext_cE_of_G(F)_by_T_J}) under $\Lambda\to G(F)$, $\lambda\mapsto t^{\lambda}$. The commutator in (\ref{ext_V_EE_of_Lambda}) is given by $(\lambda_1,\lambda_2)_c=(-1)^{\bar\kappa(\lambda_1, \lambda_2)}$, where 
$$
\bar\kappa=-\beta-\sum_{j\in J} c_j\kappa_j
$$
Let $\Gm$ act on $V_{\EE}$ via the homomoprhism $\Gm\to\Gm$, $z\mapsto z^N$. Let $\bar V_{\EE}$ be the stack quotient of $V_{\EE}$ by this action of $\Gm$. It fits into an extension of group stacks
\begin{equation}
\label{ext_bar_V_EE_of_Lambda}
1\to B(\mu_N)\to \bar V_{\EE}\to\Lambda\to 1
\end{equation}
Set
$$
\Lambda^{\sharp}=\{\lambda\in\Lambda\mid \bar\kappa(\lambda)\in N\check{\Lambda}\}
$$

We further assume that (\ref{ext_V_EE_of_Lambda}) is the push-out of the extension
\begin{equation}
\label{ext_V_EE2_of_Lambda}
1\to \mu_2\to V_{\EE, 2}\to \Lambda\to 1
\end{equation}
Recall that the exact sequence 
\begin{equation}
\label{ext_mu_2_by_mu_bar_N} 
 1\to \mu_{N}\to \mu_{2N}\to \mu_2\to 1
\end{equation} 
yields a morphism of abelian group stacks $\mu_2\to B(\mu_N)$, and the push-out of (\ref{ext_V_EE2_of_Lambda}) by this map identifies canonically with (\ref{ext_bar_V_EE_of_Lambda}). For $N$ odd the sequence (\ref{ext_mu_2_by_mu_bar_N}) splits canonically, so we get a morphism of group stacks 
\begin{equation}
\label{map_not_named_first}
\Lambda\to \bar V_{\EE}, 
\end{equation}
which is a section of (\ref{ext_bar_V_EE_of_Lambda}). Our additional input datum is a morphism for any $N$ of group stacks $\gt_{\EE}: \Lambda^{\sharp}\to \bar V_{\EE}$ extending $\Lambda^{\sharp}\hook{}\Lambda$. For $N$ odd $\gt_{\EE}$ is required to coincide with the restriction of (\ref{map_not_named_first}). For $N$ even such $\gt_{\EE}$ exists, because the restriction of (\ref{ext_V_EE_of_Lambda}) to $\Lambda^{\sharp}$ is abelian in that case.

  
\subsection{Category $\Perv_{G,\zeta}$} 

\subsubsection{} 
\label{Section_action_of_AutO}
Let $\Aut(\cO)$ be the group ind-scheme over $k$ such that , for a $k$-algebra $B$, $\Aut(\cO)(B)$ is the automorphism group of the topological $B$-algebra $B\hat\otimes \cO$ (as in \cite{FL}, Section~2.1). Let $\Aut^0(\cO)$ be the reduced part of $\Aut(\cO)$. The group scheme $\Aut^0(\cO)$ acts naturally on the exact sequence (\ref{ext_E^a_j}) acting trivially on $\Gm$ and preserving $G(\cO)$. The group scheme $\Aut^0(\cO)$ acts naturally on $F$, and the tame symbol (\ref{tame_symbol}) is $\Aut^0(\cO)$-invariant. So, by functoriality, $\Aut^0(\cO)$ acts on (\ref{ext_E_beta_of_G_ab}) acting trivially on $\Gm$. By functoriality, this gives an action of $\Aut^0(\cO)$ on (\ref{ext_cE_of_G(F)_by_T_J}) such that $\Aut^0(\cO)$ acts trivially on $\Gm$. 
 
\subsubsection{} 
\label{Section_152}
For $\lambda\in \Lambda$ let $t^{\lambda}\in \Gr_G$ denote the image of $t$ under $\lambda: F^*\to G(F)$. The set of $G(\cO)$-orbits on $\Gr_G$ idetifies with the set $\Lambda^+$ of dominant coweights of $G$. For $\lambda\in\Lambda^+$ write $\Gr^{\lambda}$ for the $G(\cO)$-orbit on $\Gr_G$ through $t^{\lambda}$. The $G$-orbit through $t^{\lambda}$ identifies with the partial flag variety $\cB^{\lambda}=G/P^{\lambda}$, where $P^{\lambda}$ is a paraboic subgroup whose Levi has the Weyl group $W^{\lambda}\subset W$ coinciding with the stabilizor of $\lambda$ in $W$. For $\lambda\in\Lambda^+$ let $\Gra^{\lambda}$ be the preimage of $\Gr^{\lambda}$ in $\Gra_G$. 
 
 The action of the loop rotation group $\Gm\subset \Aut^0(\cO)$ contracts $\Gr^{\lambda}$ to $\cB^{\lambda}\subset \Gr^{\lambda}$, we denote by $\tilde\omega_{\lambda}: \Gr^{\lambda}\to\cB^{\lambda}$ the corresponding map.  
 
 For a free $\cO$-module $\cE$ write $\cE_{\bar c}$ for its geometric fibre. Let $\Omega$ be the completed module of relative differentials of $\cO$ over $k$. For a root $\check{\alpha}$ let $\gg^{\check{\alpha}}\subset \gg$ denote the corresponding root subspace. Fix a collection $\Phi$ of trivializations $\phi_{\check{\alpha}}: \gg^{\check{\alpha}}\;\iso\; k$ for all positive roots. They yield such trivializations for all the roots of $G$ via the corresponding morphisms $\slLie_2\to \Lie(G)$.
 
  If $\check{\nu}\in\check{\Lambda}$ is orthogonal to all coroots $\alpha$ of $G$ satisfying $\<\check{\alpha},\lambda\>=0$ then we denote by $\cO(\check{\nu})$ the $G$-equivariant line bundle on $\cB^{\lambda}$ corresponding to the character  $\check{\nu}: P^{\lambda}\to\Gm$. The line bundle $\cO(\check{\nu})$ is trivialized at $1\in \cB^{\lambda}$. 
  
  Sometimes, we view $\beta$ as $\beta: \Lambda\to\check{\Lambda}$, similarly for $\kappa_j: \Lambda\to\check{\Lambda}$. The group $\Aut^0(\cO)$ acts on $\Omega_{\bar c}$ by the character denoted $\check{\epsilon}$. 
\begin{Lm} 
\label{Lm_stack_of_cL_j}
Let $\lambda\in \Lambda^+$. \\
i) For each $j\in J$ the collection $\Phi$ yields a uniquely defined $\ZZ/2\ZZ$-graded $\Aut^0(\cO)$-equivariant isomorphism 
$$
\cL_j\mid_{\Gr^{\lambda}}\,\iso\, \Omega_{\bar c}^{\frac{\kappa_j(\lambda,\lambda)}{2}}\otimes \tilde\omega_{\lambda}^*\cO(\kappa_j(\lambda))
$$
ii) The restriction of the line bundle $E_{\beta}/G(\cO)\to\Gr_G$ to $\Gr^{\lambda}$ is constant with fibre $\epsilon^{\bar\lambda}$, where $\bar\lambda\in \Lambda_{ab}$ is the image of $\lambda$. The group $G(\cO)$ acts on it by the character $G(\cO)\to G\toup{\beta(\lambda)}\Gm$, and $\Aut^0(\cO)$ acts on it by $\check{\epsilon}^{\frac{\beta(\lambda,\lambda)}{2}}$. 
\end{Lm}
\begin{Prf}
We only give the proof of the last part of ii), the rest is left to a reader. Pick a bilinear form $B: \Lambda_{ab}\otimes \Lambda_{ab}\to\ZZ$ such that $B+{^tB}=\beta$, where $^tB(\lambda_1,\lambda_2)=B(\lambda_2,\lambda_1)$ for $\lambda_i\in \Lambda_{ab}$. For this calculation we may assume $E_{\beta}=\Gm\times G_{ab}(F)$ with the product given by $(z_1, u_1)(z_2, u_2)=(z_1z_2\bar f(u_1, u_2), u_1u_2)$ for $u_i\in G_{ab}(F), z_i\in \Gm$. Here $\bar f: G_{ab}(F)\times G_{ab}(F)\to \Gm$ is the unique bimultiplicative map such that 
$$
\bar f(\lambda_1\otimes f_1, \lambda_2\otimes f_2)=(f_1, f_2)_{st}^{-B(\lambda_1,\lambda_2)}
$$ 
Let $g\in \Aut^0(\cO)$ and $b=\check{\epsilon}(g)$. Then $g$ sends $(1, t^{\bar\lambda})$ to $(1, b^{\bar\lambda}t^{\bar\lambda})\in (\bar f(t^{\bar\lambda}, b^{\bar\lambda})^{-1}, 1)(1, t^{\bar\lambda})G_{ab}(\cO)$. Finally, $\bar f(t^{\bar\lambda}, b^{\bar\lambda})=b^{-\frac{\beta(\bar\lambda,\bar\lambda)}{2}}$.
\end{Prf}


 Set $\Lambda^{\sharp,+}=\Lambda^{\sharp}\cap\Lambda^+$. For $\lambda\in\Lambda^+$ the scheme $\Gra^{\lambda}$ admits a $G(\cO)$-equivariant local system with $\Gm$-monodromy $\zeta$ if and only if $\lambda\in \Lambda^{\sharp, +}$. 
 
 By Lemma~\ref{Lm_stack_of_cL_j}, for $\lambda\in\Lambda^+$ there is a $\Aut^0(\cO)$-equivariant isomorphism between $\Gra^{\lambda}$ and the punctured (that is, with zero section removed) total space of the line bundle 
$$
\Omega_{\bar c}^{-\frac{\bar\kappa(\lambda,\lambda)}{2}}\otimes \tilde\omega_{\lambda}^*\cO(-\bar\kappa(\lambda))
$$ 
over $\Gr^{\lambda}$. Write $\Omega^{\frac{1}{2}}(\cO)$ for the groupoid of square roots of $\Omega$.  For $\cE\in \Omega^{\frac{1}{2}}(\cO)$ and $\lambda\in \Lambda^{\sharp, +}$ define the line bundle $\cL_{\lambda,\cE}$ on $\Gr^{\lambda}$ as
$$
\cL_{\lambda,\cE}=\cE_{\bar c}^{-\frac{\bar\kappa(\lambda,\lambda)}{N}}\otimes \tilde\omega_{\lambda}^*\cO(-\frac{\bar\kappa(\lambda)}{N})
$$
Let $\ocL_{\lambda,\cE}$ denote the punctured total space of 
$\cL_{\lambda,\cE}$. Let $p_{\lambda}: \ocL_{\lambda,\cE}\to\Gra^{\lambda}$ be the map over $\Gr^{\lambda}$ sending $z$ to $z^N$.   
Let $\cW^{\lambda}_{\cE}$ be the $G(\cO)$-equivariant rank one local system on $\Gra^{\lambda}$ with $\Gm$-monodromy $\zeta$ equipped with an isomorphism $p_{\lambda}^*\cW^{\lambda}_{\cE}\,\iso\,\Qlb$. Let $\cA^{\lambda}_{\cE}\in \PPerv_{G,\zeta}$ be the intermediate extension of $\cW^{\lambda}_{\cE}[\dim\Gr^{\lambda}]$ under $\Gra^{\lambda}\hook{}\Gra_G$. The perverse sheaf $\cA^{\lambda}_{\cE}$ is defined up to a  scalar automorphism (for $\lambda \in \Ker(\Lambda\to  \Lambda_{ab})$ it is defined up to a unique isomorphism). 
 
  Let $\wt\Gr^{\lambda}$ denote the restriction of the gerb $\wt\Gr_G$ to $\Gr^{\lambda}$. For $\lambda\in \Lambda^{\sharp, +}$ the map $p_{\lambda}$ yields a section 
$
s_{\lambda}: \Gr^{\lambda}\to\wt\Gr^{\lambda}
$. 
   
\begin{Lm} 
\label{Lm_parity_vanishing}
If $\lambda\in \Lambda^{\sharp, +}$ then $\cA^{\lambda}_{\cE}$ has non-trivial usual cohomology sheaves only in degrees of the same parity.
\end{Lm}
\begin{Prf}
Let $\Fl_G$ denote the affine flag variety of $G$, $q:\Fl_G\to \Gr_G$ the projection, write $\tilde q: \wt\Fl_G\to\wt\Gr_G$ for the map obtained from $q$ by the base change $\wt\Gr_G\to\Gr_G$. It suffices to prove this parity vanishing for $\tilde q^*\cA^{\lambda}_{\cE}$, this is done in \cite{Lu}.
\end{Prf}

\medskip

 Lemma~\ref{Lm_parity_vanishing} implies as in (\cite{BD1}, Proposition~5.3.3) that the category $\PPerv_{G,\zeta}$ is semisimple. 
 
\subsection{Convolution} 
\label{Section_Convolution}
Let $\tau$ be the automorphism of $\EE\times\EE$ sending $(g,h)$ to $(g, gh)$. Let $G(\cO)\times G(\cO)\times \Gm$ act on $\EE\times\EE$ so that $(\alpha,\beta, b)$ sends $(g,h)$ to $(g\beta^{-1}b^{-1}, \beta bh\alpha)$. Write $\EE\times_{G(\cO)\times \Gm} \Gra_G$ for the quotient of $\EE\times\EE$ under this free action. Then $\tau$ induces an isomorphism
$$
\bar\tau: \EE\times_{G(\cO)\times \Gm} \Gra_G\,\iso\, \Gr_G\times \Gra_G
$$
sending $(g,hG(\cO))$ to $(\bar gG(\cO), ghG(\cO))$, where $\bar g\in G(F)$ is the image of $g\in \EE$. Let $m$ be the composition of $\bar\tau$ with the projection to $\Gra_G$. Let $p_G: \EE\to \Gra_G$ be the map $h\mapsto hG(\cO)$. As in \cite{FL}, we get a diagram
$$
\Gra_G\times\Gra_G \getsup{p_G\times \id} \EE\times \Gra_G \toup{q_G} \EE\times_{G(\cO)\times \Gm} \Gra_G\toup{m} \Gra_G,
$$
where $q_G$ is the quotient map under the action of $G(\cO)\times \Gm$. 

 For $K_i\in \PPerv_{G,\zeta}$ define the convolution $K_1\ast K_2\in \D(\Gra_G)$ by 
$
K_1\ast K_2=m_! K\in\D(\Gra_G),
$
where $K[1]$ is a perverse sheaf on $\EE\times_{G(\cO)\times \Gm} \Gra_G$ equipped with an isomorphism 
$$
q_G^*K\,\iso\, p_G^*K_1\boxtimes K_2
$$ 
Since $q_G$ is a $G(\cO)\times \Gm$-torsor, and $p_G^*K_1\boxtimes K_2$ is naturally equivariant under $G(\cO)\times \Gm$-action, $K$ is defined up to a unique isomorphism. As in (\cite{FL}, Lemma~2.6), one shows that $K_1\ast K_2\in  \PPerv_{G,\zeta}$. 

 For $K_i\in \PPerv_{G,\zeta}$ one similarly defines the convolution $K_1\ast K_2\ast K_3\in \PPerv_{G,\zeta}$ and shows that $(K_1\ast K_2)\ast K_3\,\iso\, K_1\ast K_2\ast K_3\,\iso\, K_1\ast (K_2\ast K_3)$ canonically. Besides, $\cA^0_{\cE}$ is a unit object in $\PPerv_{G,\zeta}$. 
 
\subsection{Fusion} 
\label{Section_Fusion}

As in \cite{FL}, we are going to show that the convolution product on $\PPerv_{G,\zeta}$ can be interpreted as a fusion product, thus leading to a commutativity constraint on $\PPerv_{G,\zeta}$.  
 
 Fix $\cE\in\Omega^{\frac{1}{2}(\cO)}$. Let $\Aut_2(\cO)=\Aut(\cO,\cE)$ be the group scheme defined in (\cite{FL}, Section~2.3), let $\Aut_2^0(\cO)$ be the preimage of $\Aut^0$ in $\Aut_2(\cO)$. 
 
 Let $\lambda\in\Lambda^{\sharp, +}$. Since $p_{\lambda}: \ocL_{\lambda,\cE}\to \Gra^{\lambda}$ is $\Aut_2^0(\cO)$-equivariant, the action of $\Aut^0(\cO)$ on $\Gra_G$ lifts to a $\Aut_2^0(\cO)$-equivariant structure on $\cA^{\lambda}_{\cE}$. As in (\cite{FL}, Section~2.3) one shows that the corresponding $\Aut_2^0(\cO)$-equivariant structure on each $\cA^{\lambda}_{\cE}$ is unique.
 
  For $x\in X$ let $\cO_x$ be the completed local ring at $x\in X$, $F_x$ its fraction field. Write $\cF^0_G$ for the trivial $G$-torsor on a base. Write $\Gr_{G,x}=G(F_x)/G(\cO_x)$ for the corresponding affine grassmanian. Recall that $\Gr_{G,x}$ can be seen as the ind-scheme classifying a $G$-torsor $\cF$ on $X$ together with a trivialization $\nu: \cF\,\iso\, \cF^0_G\mid_{X-x}$. 
  
  For $m\ge 1$ let $\Gr_{G,X^m}$ and $G_{X^m}$ be defined as in (\cite{FL}, Section~2.3). Recall that $\Gr_{G, X^m}$ is the ind-scheme classifying $(x_1,\ldots, x_m)\in X^m$, a $G$-torsor $\cF_G$ on $X$, and a trivialization $\cF_G\,\iso\, \cF^0_G\mid_{X-\cup x_i}$. Here $G_{X^m}$ is a group scheme over $X^m$ classifying $\{(x_1,\ldots, x_m)\in X^m,\mu\}$, where $\mu$ is an automorphism of $\cF^0_G$ over the formal neighbourhood of $D=\cup_i x_i$ in $X$. 
  
  For $j\in J$ let $\cL_{j, X^m}$ be the ($\ZZ/2\ZZ$-graded purely of parity zero)  line bundle on $\Gr_{G, X^m}$ whose fibre at $(\cF_G, x_1,\ldots, x_m)$ is 
$$
\det\RG(X, (\gg_j)_{\cF^0_G})\otimes \det\RG(X, (\gg_j)_{\cF_G})^{-1}
$$ 
Here for a $G$-module $V$ and a $G$-torsor $\cF_G$ on a base $S$ we write $V_{\cF_G}$ for the induced vector bundle on $S$. 
  
  As in Section~\ref{Section_notations}, our choice of $\Omega^{\frac{1}{2}}$ yields a functor 
\begin{equation}  
\label{functor_creating_theta_datum}  
\cE^s(G_{ab})\to \cP^{\theta}(X, \Lambda_{ab})
\end{equation}  
Let $\theta_0\in \cP^{\theta}(X, \Lambda_{ab})$ denote the image under this functor of the extension (\ref{ext_Lambda_ab_by_Gm}) with the bilinear form $-\beta$.
  
  For a reductive group $H$ write $\Bun_H$ for the stack of $H$-torsors on $X$.  Write $\cPic(\Bun_H)$ for the groupoid of super line bundles on $\Bun_H$. 
For $\mu\in \pi_1(H)$ write $\Bun_H^{\mu}$ for the connected component of $\Bun_H$ classifying $H$-torsors of degree $-\mu$. Similarly, for $\mu\in \pi_1(G)$ we denote by $\Gr_G^{\mu}$ the connected component containing $t^{\lambda}G(\cO)$ for any $\lambda\in\Lambda$ over $\mu$.
  
 Recall the functor $\cP^{\theta}(X, \Lambda_{ab})\to \cPic(\Bun_{G_{ab}})$ defined in (\cite{L}, Section~4.2.1, formula (18)). Let $\cL_{\beta}\in \cPic(\Bun_{G_{ab}})$ denote the image of $\theta_0$ under this functor. It is purely of parity zero as $\ZZ/2\ZZ$-graded. For $\mu\in \Lambda_{ab}$ we have a map $i_{\mu}:X\to \Bun_{G_{ab}}$, $x\mapsto \cO(-\mu x)$. By definition, 
$$
i_{\mu}^*\cL_{\beta}\,\iso\, (\Omega^{\frac{1}{2}})^{\beta(\mu,\mu)}\otimes\epsilon^{\mu}
$$ 

 For $m\ge 1$ let $\cL_{\beta, X^m}$ be the pull-back of $\cL_{\beta}$ under $\Gr_{G, X^m}\to\Bun_{G_{ab}}$. Let  $\Gra_{G, X^m}$ denote the punctured total space of the line bundle over $\Gr_G$
$$
\cL_{\beta, X^m}\otimes(\mathop{\otimes}\limits_{j\in J} (\cL_{j, X^m})^{c_j})
$$

\begin{Rem} The line bundle $\cL_{\beta, X^m}$ is $G_{X^m}$-equivariant. For $(x_1,\ldots, x_m)\in X^m$ let 
$$
\{y_1,\ldots, y_s\}=\{x_1,\ldots, x_m\}
$$ 
with $y_i$ pairwise different. Let $\mu_i\in\Lambda$ for $1\le i\le s$. Consider a point $\eta\in \Gr_{G, X^m}$ over $\bar\eta\in \Gr_{G_{ab}, X^m}$ given by $\cF^0_{G_{ab}}(-\sum_{i=1}^s\mu_iy_i)$ with the evident trivialization over $X-\cup_i y_i$. The fibre of $G_{X^m}$ at $(x_1,\ldots, x_m)$ is $\prod_{i=1}^s G(\cO_{y_i})$, this group acts on the fibre $(\cL_{\beta, X^m})_{\eta}$ by the character
$$
\prod_{i=1}^s G(\cO_{y_i})\to \prod_{i=1}^sG_{ab}(\cO_{y_i})\to \prod_{i=1}^s G_{ab}\stackrel{\prod_i\beta(\mu_i)}{\longrightarrow} \Gm
$$ 
Since the line bundles $\cL_{j, X^m}$ are also $G_{X^m}$-equivariant, the action of $G_{X^m}$ on $\Gr_{G, X^m}$ is lifted to an action on $\Gra_{G, X^m}$. 
\end{Rem}

 Let $\Perv_{G, \zeta, X^m}$ be the category of $G_{X^m}$-equivariant perverse sheaves on $\Gra_{G, X^m}$ with $\Gm$-monodromy $\zeta$. Set
$$
\PPerv_{G,\zeta, X^m}=\Perv_{G, \zeta, X^m}[-m-1]\subset \D(\Gra_{G, X^m})
$$

  For $x\in X$ let $D_x=\Spec\cO_x$, $D_x^*=\Spec F_x$. The analog of the convolution diagram from (\cite{FL}, Section~2.3) is the following one, where the left and right squares are cartesian:
$$
\begin{array}{ccccccc}
\Gra_{G,X}\times \Gra_{G,X} & \getsup{\tilde p_{G,X}} & \tilde C_{G,X} & \toup{\tilde q_{G,X}} & \wt\Conv_{G,X} & \toup{\tilde m_X} & \Gra_{G, X^2}\\
\downarrow &&  \downarrow &&  \downarrow &&  \downarrow \\
\Gr_{G,X}\times \Gr_{G,X} & \getsup{p_{G,X}} & C_{G,X} & \toup{q_{G,X}} & \Conv_{G,X} & \toup{m_X} & \Gra_{G, X^2}
\end{array}
$$  
Here the low row is the convolution diagram from \cite{FL}. Namely, $C_{G,X}$ is the ind-scheme classifying collections:
\begin{equation}
\label{collection_C_GX}
x_1,x_2\in X, \; G\mbox{-torsors}\; \cF^1_G,\cF^2_G\; \mbox{on} \; X\; \mbox{with}\; \nu_i: \cF^i_G\,\iso\, \cF^0_G\mid_{X-x_i}, \; \mu_1: \cF^1_G\,\iso\, \cF^0_G\mid_{D_{x_2}}
\end{equation}
The map $p_{G,X}$ forgets $\mu_1$. The ind-scheme $\Conv_{G,X}$ classifies collections:
\begin{multline}
\label{collection_Conv_GX}
x_1,x_2\in X, \; G\mbox{-torsors}\; \cF^1_G,\cF^2_G\; \mbox{on} \; X,\\  
\mbox{isomorphisms}\; \nu_1: \cF^1_G\,\iso\, \cF^0_G\mid_{X-x_1}\; \mbox{and}\; \eta: \cF^1_G\,\iso\, \cF_G\mid_{X-x_2}
\end{multline}
The map $m_X$ sends this collection to $(x_1, x_2, \cF_G)$ together with the trivialization $\eta\comp\nu_1^{-1}: \cF^0_G\,\iso\, \cF_G\mid_{X-x_1-x_2}$. 
 
 The map $q_{G,X}$ sends (\ref{collection_C_GX}) to (\ref{collection_Conv_GX}), where $\cF_G$ is obtained by gluing $\cF^1_G$ on $X-x_2$ and $\cF^2_G$ on $D_{x_2}$ using their identification over $D_{x_2}^*$ via $\nu_2^{-1}\comp\mu_1$. 
  
  For $j\in J$ there is a canonical $\ZZ/2\ZZ$-graded isomorphism
\begin{equation}
\label{iso_for_j_factorization}
q_{G,X}^*m_X^*\cL_{j, X^2}\,\iso\, p_{G,X}^*(\cL_{j, X}\boxtimes \cL_{j, X})
\end{equation}  
\begin{Lm} There is a canonical $\ZZ/2\ZZ$-graded isomorphism
\begin{equation}
\label{iso_for_beta_factorization}
q_{G,X}^*m_X^*\cL_{\beta, X^2}\,\iso\, p_{G,X}^*(\cL_{\beta, X}\boxtimes \cL_{\beta, X})
\end{equation}
\end{Lm}
\begin{Prf}
This isomorphism comes from the corresponding isomorphism for $G_{ab}$, so for this proof we may assume $G=G_{ab}$. For a point (\ref{collection_C_GX}) of $C_{G,X}$ consider its image under $q_{G,X}$ given by (\ref{collection_Conv_GX}). Note that $\cF_G=\cF^1_G\otimes\cF^2_G$ with the trivialization $\nu_1\otimes\nu_2: \cF^1_G\otimes\cF^2_G\,\iso\,\cF^0_G\mid_{X-x_1-x_2}$. One gets by (\cite{L}, Proposition~4.2)
$$
(\cL_{\beta})_{(\cF^1_G, \nu_1)}\otimes (\cL_{\beta})_{(\cF^2_G, \nu_2)}\,\iso\, (\cL_{\beta})_{\cF^1_G}\otimes (\cL_{\beta})_{\cF^2_G} \otimes (^{-\beta}\cL^{univ}_{\cF^1_G, \cF^2_G})\,\iso\, (\cL_{\beta})_{\cF^1_G\otimes\cF^2_G}
$$
with the notations of \select{loc.cit.} Here we used the following trivialization $(^{-\beta}\cL^{univ}_{\cF^1_G, \cF^2_G})\,\iso\, k$. Forgetting about nilpotents for simplicity, we may assume $\cF^2_G\,\iso\, \cF^0_G(\nu x_2)$ for some $\nu\in\Lambda$ with the evident trivialization over $X-x_2$. Then 
$$
(^{-\beta}\cL^{univ}_{\cF^1_G, \cF^2_G})\,\iso\, (\cL^{-\beta(\nu)}_{\cF^1_G})_{x_2}\,\iso\, k,
$$ 
the latter isomorphism is obtained from $\mu_1: \cF^1_G\,\iso\, \cF^0_G\mid_{D_{x_2}}$.
\end{Prf}

\medskip

 The isomorphisms (\ref{iso_for_j_factorization}) and (\ref{iso_for_beta_factorization}) allow to define the map $\tilde q_{G,X}$ exactly as in (\cite{FL}, Section~2.3), this is the product of the corresponding maps. 
 
  Now for $K_i\in \PPerv_{G, \zeta, X}$ there is a (defined up to a unique isomorphism) perverse sheaf $K_{12}[3]$ on $\wt\Conv_{G,X}$ equipped with
$\tilde q_{G,X}^*K_{12}\,\iso\, \tilde p_{G,X}^*(K_1\boxtimes K_2)$. Moreover $K_{12}$ has $\Gm$-monodromy $\zeta$. We let
$$
K_1\ast_X K_2=\tilde m_{X !} K_{12}
$$   
As in (\cite{FL}, Section~2.3) one shows that $K_1\ast_X K_2\in \PPerv_{G, \zeta, X^2}$. 

 Let $\cE\in\Omega^{\frac{1}{2}}(\cO)$. As in \select{loc.cit.}, one has the $\Aut^0_2(\cO)$-torsor $\hat X_2\to X$ whose fibre over $x$ is the scheme of isomorphisms between $(\Omega^{\frac{1}{2}}_x, \cO_x)$ and $(\cE, \cO)$. One has the isomorphisms 
$$
\Gr_{G,X}\,\iso\, \hat X_2\times_{\Aut^0_2(\cO)} \Gr_G\;\;\mbox{and}\;\; \Gra_{G,X}\,\iso\, \hat X_2\times_{\Aut^0_2(\cO)} \Gra_G
$$
Since any $K\in \Perv_{G,\zeta}$ is $\Aut^0_2(\cO)$-equivariant, we get the fully faithful functor 
$$
\tau^0: \PPerv_{G,\zeta}\to\PPerv_{G,\zeta, X}
$$ 
sending $K$ to the descent of $\Qlb\boxtimes K$ under $\hat X_2\times\Gra_G\to \Gra_{G,X}$.

 Let $U\subset X^2$ be the complement to the diagonal. Let $j: \Gra_{G,X^2}(U)\hook{} \Gra_{G,X^2}$ be the preimage of $U$. Let $i: \Gra_{G,X}\to\Gra_{G,X^2}$ be obtained by the base change $X\to X^2$. Recall that $\tilde m_X$ is an isomorphism over $\Gra_{G,X^2}(U)$. For $F_i\in \PPerv_{G,\zeta}$ letting $K_i=\tau^0 F_i$ define
$$
K_{12}\mid_U:=K_{12}\mid_{\Gra_{G,X^2}(U)}
$$
as above, it is placed in perverse degree $3$. Then $K_1\ast_X K_2\,\iso\, j_{!*}(K_{12}\mid_U)$ and $\tau^0(F_1\ast F_2)\,\iso\, i^*(K_1\ast_X K_2)$. So, the involution $\sigma$ of $\Gra_{G,X^2}$ interchanging $x_i$ yields 
$$
\tau^0(F_1\ast F_2)\,\iso\, i^* j_{!*}(K_{12}\mid_U)\,\iso\, i^* j_{!*}(K_{21}\mid_U)\,\iso\, \tau^0(F_2\ast F_1),
$$
because $\sigma^*(K_{12}\mid_U)\,\iso\, K_{21}\mid_U$ canonically. As in \cite{FL}, the associativity and commutativity constraints are compatible, so $\PPerv_{G,\zeta}$ is a symmetric monoidal category. 

\begin{Rem} 
\label{Rem_rigid}
Let $\P_{G(\cO)}(\Gra_G)$ denote the category of $G(\cO)$-equivariant perverse sheaves on $\Gra_G$. One has the covariant self-functor $\star$ on $\P_{G(\cO)}(\Gra_G)$ induced by the map $\EE\to\EE$, $z\mapsto z^{-1}$. Then $K\mapsto K^{\vee}:=\DD(\star K)[-2]$ is a contravariant functor $\PPerv_{G,\zeta}\to \PPerv_{G,\zeta}$. As in (\cite{FL}, Remark~2.8), one shows that $\RHom(K_1\ast K_2, K_3)\,\iso\, \RHom(K_1, K_3\ast K_2^{\vee})$. So, $K_3\ast K_2^{\vee}$ represents the internal $\HOM(K_2, K_3)$ in the sense of the tensor structure on $\PPerv_{G,\zeta}$. Besides, $\star(K_1\ast K_2)\,\iso\, (\star K_2)\ast (\star K_1)$ canonically.
\end{Rem}

\subsection{Main result}
\label{Section_MR}

Below we introduce a tensor category $\PPerv^{\natural}_{G,\zeta}$ obtained from $\PPerv_{G,\zeta}$ by some modification of the commutativity constraint. Let $\check{T}_{\zeta}=\Spec k[\Lambda^{\sharp}]$ be the torus whose weight lattice is $\Lambda^{\sharp}$. 

 For $a\in \QQ^*$ written as $a=a_1/a_2$ with $a_i\in\ZZ$ prime to each other and $a_2>0$, say that $a_2$ is the denominator of $a$. Recall that we assume $N$ invertible in $k$.

\begin{Th} 
\label{Th_one}
There is a connected reductive group $\check{G}_{\zeta}$ over $\Qlb$ and a canonical equivalence of tensor categories 
$$
\PPerv^{\natural}_{G,\zeta}\,\iso\, \Rep(\check{G}_{\zeta}) .
$$
There is a canonical inclusion $\check{T}_{\zeta}\subset \check{G}_{\zeta}$ whose image is a maximal torus in $\check{G}_{\zeta}$. The Weyl groups of $G$ and $\check{G}_{\zeta}$ viewed as subgroups of $\Aut(\Lambda^{\sharp})$ are the same. Our choice of a Borel subgroup $T\subset B\subset G$ yields a Borel subgroup $\check{T}_{\zeta}\subset \check{B}_{\zeta}\subset \check{G}_{\zeta}$. 
The corresponding simple roots (resp., coroots) of $(\check{G}_{\zeta}, \check{T}_{\zeta})$ are $\delta_i\alpha_i$ (resp., $\check{\alpha}_i/\delta_i$) for $i\in\cI$. Here $\delta_i$ is the denominator of $\frac{\bar\kappa(\alpha_i,\alpha_i)}{2N}$. 
\end{Th}

\begin{Rem} i) The root datum described in Theorem~\ref{Th_one} is defined uniquely. The roots are the union of $W$-orbits of simple roots. For $\alpha\in R$ let $\delta_{\alpha}$ denote the denominator of $\frac{\bar\kappa(\alpha,\alpha)}{2N}$. Then $\delta_{\alpha}\alpha$ is a root of $\check{G}_{\zeta}$. Any root of $\check{G}_{\zeta}$ is of this form. Compare with the metaplectic root datum appeared in (\cite{Mc}, \cite{W}, \cite{R}).\\
ii) We hope there could exist an improved construction, which is a functor from the category of central extensions $1\to K_2\to E\to G\to 1$ over $k$ to the 2-category of symmetric monoidal categories, $E\mapsto \PPerv_{G, E}$ such that $\PPerv_{G,E}$ is tensor equivalent to the category $\Rep(\check{G}_E)$ of representations of some connected reductive group $E$. 
\end{Rem}

\section{Proof of Theorem~\ref{Th_one}}

\subsection{Functors $F'_P$} 
\label{Section_functors_F'_P}
Let $P\subset G$ be a parabolic subgroup containing $B$. Let $M\subset P$ be its Levi factor containing $T$. Let $\cI_M\subset \cI$ be the subset parametrizing the simple roots of $M$. Write
$$
1\to \Gm\to \EE_M\to M(F)\to 1
$$
for the restriction of (\ref{ext_cE_of_G(F)_by_T_J}) to $M(F)$. It is equipped with an action of $\Aut^0(\cO)$ and a section over $M(\cO)$ coming from the corresponding objects for (\ref{ext_cE_of_G(F)_by_T_J}). 

 Write $\Gr_M, \Gr_P$ for the affine grassmanians for $M, P$ respectively. For $\theta\in \pi_1(M)$ write $\Gr_M^{\theta}$ for the connected component of $\Gr_M$ containing $t^{\lambda}M(\cO)$ for any $\lambda\in\Lambda$ over $\theta\in\pi_1(M)$. The diagram $M\gets P\to G$ yields the following diagram of affine grassmanians
$$
\Gr_M \getsup{\gt_P}\Gr_P\toup{\gs_P} \Gr_G \, .
$$
Let $\Gr_P^{\theta}$ be the connected component of $\Gr_P$ such that $\gt_P$ restricts to a map $\gt^{\theta}_P: \Gr^{\theta}_P\to \Gr^{\theta}_M$. Write $\gs^{\theta}_P: \Gr^{\theta}_P\to \Gr_G$ for the restriction of $\gs_P$. The restriction of $\gs^{\theta}_P$ to $(\Gr^{\theta}_P)_{red}$ is a locally closed immesion. 

 The section $M\to P$ yields a section $\gr_P: \Gr_M\to \Gr_P$ of $\gt_P$. By abuse of notations, write
$$
\Gra_M\toup{\gr_P} \Gra_P\toup{\gs_P} \Gra_G
$$
for the diagram obtained from $\Gr_M\toup{\gr_P} \Gr_P\toup{\gs_P}\Gr_G$ by the base change $\Gra_G\to \Gr_G$. Note that $\gt_P$ lifts naturally to a map denoted $\gt_P: \Gra_P\to \Gra_M$ by abuse of notations.  

 Let $\Perv_{M, G, \zeta}$ denote the category of $M(\cO)$-equivariant perverse sheaves on $\Gra_M$ with $\Gm$-monodromy $\zeta$. Set
$$
\PPerv_{M,G,\zeta}=\Perv_{M,G,\zeta}[-1]\subset \D(\Gra_M) \, .
$$
Define the functor
$$
F'_P: \PPerv_{G,\zeta}\to \D(\Gra_M)
$$
by $F'_P(K)=\gt_{P !}\gs^*_PK$. Write $\Gra_M^{\theta}$ for the connected component of $\Gra_M$ over $\Gr_M^{\theta}$, similarly for $\Gra_P^{\theta}$. 
Write 
$$
\PPerv^{\theta}_{M,G, \zeta}\subset \PPerv_{M,G,\zeta}
$$
for the full subcategory of objects that vanish off $\Gra^{\theta}_M$. Set
$$
\PPerv'_{M,G, \zeta}=\mathop{\oplus}_{\theta\in \pi_1(M)} \PPerv^{\theta}_{M,G,\zeta}[\<\theta, 2\check{\rho}_M-2\check{\rho}\>] \, .
$$
As in \cite{FL}, one shows that $F'_P$ sends $\PPerv_{G,\zeta}$ to $\PPerv'_{M,G, \zeta}$. This is a combination of the hyperbolic localization argument (\cite{MV}, Theorem~3.5) or (\cite{L4}, Proposition~12) with the dimension estimates of (\cite{MV}, Theorem~3.2) or (\cite{BG}, Proposition~4.3.3). 

 For the Borel subgroup $B$ the above construction gives $F'_B: \PPerv_{G,\zeta}\to \PPerv'_{T, G,\zeta}$. 
 
 Let $B(M)\subset M$ be a Borel subgroup such that the preimage of $B(M)$ under $P\to M$ equals $B$. The inclusions $T\subset B(M)\subset M$ yield a diagram
\begin{equation}
\label{diag_Gr_for_M} 
\Gr_T \toup{\gr_{B(M)}} \Gr_{B(M)} \toup{\gs_{B(M)}} \Gr_M \, .
\end{equation}
Write
$$
\Gra_T \toup{\gr_{B(M)}} \Gra_{B(M)} \toup{\gs_{B(M)}} \Gra_M
$$
for the diagram obtained from (\ref{diag_Gr_for_M}) by the base change $\Gra_M\to \Gr_M$. The projection $B(M)\to T$ yields $\gt_{B(M)}: \Gr_{B(M)}\to \Gr_T$, it lifts naturally to the map denoted $\gt_{B(M)}: \Gra_{B(M)}\to \Gra_T$ by abuse of notations. For $K\in \PPerv'_{M,G,\zeta}$ set
$$
F'_{B(M)}(K)=(\gt_{B(M)})_!\gs^*_{B(M)}K \, .
$$
As in \cite{FL}, this defines the functor $F'_{B(M)}: \PPerv'_{M,G,\zeta}\to \PPerv'_{T, G,\zeta}$, and one has canonically
\begin{equation}
\label{iso_from_P_to_B}
F'_{B(M)}\comp F'_P\,\iso\, F'_B  \, .
\end{equation}
 
\subsubsection{} 
\label{Section_2.1.1} For $j\in J$ let $\cL_{j,M}$ denote the restriction of $\cL_j$ under $\gs_P\gr_P:\Gr_M\to \Gr_G$. Let $\Lambda^+_M$ denote the coweights dominant for $M$. For $\lambda\in\Lambda^+_M$ denote by $\Gr^{\lambda}_M$ the $M(\cO)$-orbit through $t^{\lambda}M(\cO)$. Let $\Gra_M^{\lambda}$ be the preimage of $\Gr_M^{\lambda}$ under $\Gra_M\to \Gr_M$. The $M$-orbit through $t^{\lambda}M(\cO)$ is isomorphic to the partial flag variety $\cB^{\lambda}_M=M/P^{\lambda}_M$, where the Levi subgroup of $P^{\lambda}_M$ has the Weyl group coinciding with the stabilizer of $\lambda$ in $W_M$. Here $W_M$ is the Weyl group of $M$. As for $G$, we have a natural map $\tilde\omega_{M,\lambda}: \Gr^{\lambda}_M\to \cB^{\lambda}_M$. 

 If $\check{\nu}\in\check{\Lambda}$ is orthogonal to all coroots $\alpha$ of $M$ satisfying $\<\check{\alpha},\lambda\>=0$ then we denote by $\cO(\check{\nu})$ the $M$-equivariant line bundle on $\cB^{\lambda}_M$ corresponding to the character $\check{\nu}: P^{\lambda}_M\to\Gm$.  As in Lemma~\ref{Lm_stack_of_cL_j}, for $j\in J$ the collection $\Phi$ yields a uniquely defined $\ZZ/2\ZZ$-graded $\Aut^0(\cO)$-equivariant isomorphism
$$
\cL_{j,M}\mid_{\Gr^{\lambda}_M}\,\iso\, \Omega_{\bar c}^{\frac{\kappa_j(\lambda,\lambda)}{2}}\otimes\tilde\omega_{M,\lambda}^*\cO(\kappa_j(\lambda))
$$
So, for $\lambda\in\Lambda^+_M$ there is a $\Aut^0(\cO)$-equivariant isomorphism between $\Gra^{\lambda}_M$ and the punctured total space of the line bundle
$$
\Omega_{\bar c}^{-\frac{\bar\kappa(\lambda,\lambda)}{2}}\otimes\tilde\omega_{M,\lambda}^*\cO(-\bar\kappa(\lambda))
$$
over $\Gr_M^{\lambda}$. Set $\Lambda^{\sharp,+}_M=\Lambda^{\sharp}\cap \Lambda^+_M$. As for $G$ itself, for $\lambda\in \Lambda^+_M$ the scheme $\Gra_M^{\lambda}$ admits a $M(\cO)$-equivariant local system with $\Gm$-monodromy $\zeta$ if and only if $\lambda\in \Lambda^{\sharp, +}_M$. 

 As in Section~\ref{Section_152} pick $\cE\in \Omega^{\frac{1}{2}}(\cO)$. For $\lambda\in \Lambda^{\sharp, +}_M$ define the line bundle $\cL_{\lambda, M,\cE}$ on $\Gr^{\lambda}_M$ as
$$
\cL_{\lambda, M,\cE}=\cE_{\bar c}^{-\frac{\bar\kappa(\lambda,\lambda)}{N}}\otimes \tilde\omega_{M,\lambda}^*\cO(-\frac{\bar\kappa(\lambda)}{N})
$$
Let $\ocL_{\lambda, M,\cE}$ be the punctured total space of $\cL_{\lambda, M,\cE}$. Let $p_{\lambda, M}: \ocL_{\lambda, M,\cE}\to \Gra^{\lambda}_M$  be the map over $\Gr^{\lambda}_M$ sending $z$ to $z^N$. Let $\cW^{\lambda}_{M,\cE}$ be the rank one $M(\cO)$-equivariant local system $\cW^{\lambda}_{M,\cE}$ on $\Gra^{\lambda}_M$ with $\Gm$-monodromy $\zeta$ equipped with an isomorphism $p_{\lambda, M}^*\cW^{\lambda}_{M,\cE}\,\iso\, \Qlb$. Let $\cA^{\lambda}_{M, \cE}\in \PPerv_{M,G,\zeta}$ be the intermediate extension of $\cW^{\lambda}_{M,\cE}[\dim\Gr^{\lambda}_M]$ to $\Gra_M$, it is defined up to a scalar automorphism. 

Set 
$$
\wt\Gr_M=\Gr_M\times_{\Gr_G} \wt\Gr_G
$$ 
For $\lambda\in \Lambda^+_M$ let $\wt\Gr^{\lambda}_M$ be the restriction of the gerb $\wt\Gr_M$ to $\Gr^{\lambda}_M$. As for $G$ itself, for $\lambda\in \Lambda^{\sharp, +}_M$ the map $p_{\lambda, M}$ yields a section $s_{\lambda, M}: \Gr^{\lambda}_M\to \wt\Gr^{\lambda}_M$. 
 
  The analog of Lemma~\ref{Lm_parity_vanishing} holds for the same reasons. The perverse sheaf $\cA^{\lambda}_{M,\cE}$ has non-trivial cohomology sheaves only in degrees of the same parity. It follows that $\PPerv_{M,G,\zeta}$ is semisimple.
 
\subsubsection{More tensor structures} One equips $\PPerv_{M,G,\zeta}$ and $\PPerv'_{M,G,\zeta}$ with a convolution product as in Section~\ref{Section_Convolution}. The convolution for these categories can be interpreted as fusion, and this allows to define a commutativity constraint on these categories via fusion. 

 Each of the line bundles $\cL_{j, X^m}, \cL_{\beta, X^m}$ on $\Gr_{G, X^m}$ admits the factorization structure as in (\cite{FL}, Section~4.1.2). 
 
 As for $G$, we have the ind-scheme $\Gr_{M, X^m}$ for $m\ge 1$ and the group scheme $M_{X^m}$ over $X^m$ defined similarly. Let $\Gra_{M,G, X^m}\to \Gra_{G, X^m}$ be obtained from $\Gr_{M, X^m}\to \Gr_{G, X^m}$ by the base change $\Gra_{G, X^m}\to \Gr_{G, X^m}$. The group scheme $M_{X^m}$ acts naturally on $\Gra_{M, G, X^m}$. 
 
 Write $\Perv_{M, G, \zeta, X^m}$ be the category of $M_{X^m}$-equivariant perverse sheaves on $\Gra_{M,G, X^m}$ with $\Gm$-monodromy $\zeta$. Set
$$
\PPerv_{M, G, \zeta, X^m}=\Perv_{M, G, \zeta, X^m}[-m-1]\, .
$$
  
  Let $\Aut^0_2(\cO)$ act on $\Gra_M$ via its quotient $\Aut^0(\cO)$. Then every object of $\Perv_{M, G, \zeta}$ admits a unique $\Aut^0_2(\cO)$-equivariant structure. On has
$$
\Gra_{M, G, X}\,\iso\, \hat X_2\times_{\Aut^0_2(\cO)}\Gra_M,
$$
and as above one gets a fully faithful functor 
$$
\tau^0: \PPerv_{M,G, \zeta}\to \PPerv_{M,G,\zeta, X}
$$
Define the commutativity constraint on $\PPerv_{M,G,\zeta}$ and $\PPerv'_{M,G,\zeta}$ via fusion as in Section~\ref{Section_Fusion}. As in \cite{FL}, one checks that $\PPerv_{M,G,\zeta}$ and $\PPerv'_{M,G,\zeta}$ are symmetric monoidal categories. Exactly as in (\cite{FL}, Lemma~4.1), one proves the following.  

\begin{Lm} The functors $F'_P$,  $F'_{B(M)}, F'_B$ are tensor functors, and (\ref{iso_from_P_to_B}) is an isomorphism of tensor functors. \QED
\end{Lm}

\subsection{Fiber functor} Recall from Section~\ref{Section_2.1.1} that for $\lambda\in \Lambda$ the scheme $\Gra_T^{\lambda}$ admits a $T(\cO)$-equivariant local system with $\Gm$-monodromy $\zeta$ if and only if $\lambda\in \Lambda^{\sharp}$. View an object of $\PPerv_{T,G,\zeta}$ as a complex on $\wt\Gr_T$. The map $\gt_{\EE}$ from Section~\ref{Section_Input_data} defines for each $\lambda\in\Lambda^{\sharp}$ a section $\gt_{\lambda,\EE}: \Gr^{\lambda}_T\to \wt\Gr^{\lambda}_T$. For $K\in \PPerv_{T,G,\zeta}$ the complex $\gt_{\lambda,\EE}^*K$ is constant and placed in degree zero, so we view as a vector space denoted $F^{\lambda}_T(K)$. Let 
$$
F_T=\mathop{\oplus}_{\lambda\in \Lambda^{\sharp}} F^{\lambda}_T: \PPerv_{T,G,\zeta}\to\Vect
$$ 
This is a fibre functor on $\PPerv_{T,G,\zeta}$. By (\cite{DM}, Theorem~2.11) we get
$$
\PPerv_{T,G,\zeta}\,\iso\, \Rep(\check{T}_{\zeta})\, .
$$

 For $\nu\in \Lambda^{\sharp}$ write $F'^{\nu}_{B(M)}$ for the functor $F'_{B(M)}$ followed by restriction to $\Gra^{\nu}_T$. Write $F^{\nu}_M: \PPerv_{M,G,\zeta}\to \Vect$ for the functor
$$
F^{\nu}_T F'^{\nu}_{B(M)}[\<\nu, 2\check{\rho}_M\>]
$$
In particular, this definition applies for $M=G$ and gives the functor $F^{\nu}_G: \PPerv_{G, \zeta}\to \Vect$. 

 For $\nu\in \Lambda$ as in Section~\ref{Section_functors_F'_P} one has the map $\gt^{\nu}_{B(M)}: \Gr^{\nu}_{B(M)}\to \Gr^{\nu}_T$. Let $\wt\Gr^{\nu}_{B(M)}$ denote the restriction of the gerb $\wt\Gr_M$ under $\Gr^{\nu}_{B(M)}\to \Gr_M$. For $\nu\in \Lambda^{\sharp}$ the section $\gt_{\nu,\EE}: \Gr_T^{\nu}\to \wt\Gr_T^{\nu}$ yields by restriction under $\gt^{\nu}_{B(M)}$ the section that we denote $\gt_{\nu, B(M)}: \Gr^{\nu}_{B(M)}\to \wt\Gr^{\nu}_{B(M)}$.  
 
\begin{Lm} 
\label{Lm_instead_of_Lm4.2}
If $\nu\in \Lambda^{\sharp}$, $\lambda\in \Lambda^{\sharp, +}_M$ then $F^{\nu}_M(\cA^{\lambda}_{M,\cE})$ has a canonical base consisting of those irreducible components of
$$
\Gr^{\nu}_{B(M)}\cap \Gr^{\lambda}_M
$$
over which the (shifted) local system $\gt_{\nu, B(M)}^*\cA^{\lambda}_{M,\cE}$ is constant. Here we view $\cA^{\lambda}_{M,\cE}$ as a perverse sheaf on $\wt\Gr_M$. In particular,  for $w\in W_M$ one has 
$$
F^{w(\lambda)}_M(\cA^{\lambda}_{M,\cE})\,\iso\, \Qlb
$$
\end{Lm}
\begin{Prf} 
Exactly as in (\cite{FL}, Lemma~4.2). 
\end{Prf} 

\medskip\smallskip

 Consider the following $\ZZ/2\ZZ$-grading on $\PPerv'_{M, G,\zeta}$. For $\theta\in \pi_1(M)$ call an object of $\PPerv^{\theta}_{M,G,\zeta}[\<\theta, 2\check{\rho}_M-2\check{\rho}\>]$ of parity $\<\theta, 2\check{\rho}\>\mod 2$, the latter expression depends only on the image of $\theta$ in $\pi_1(G)$. As in \cite{FL}, this $\ZZ/2\ZZ$-grading on $\PPerv'_{M, G,\zeta}$ is compatible with the tensor structure. In particular, for $M=G$ we get a $\ZZ/2\ZZ$-grading on $\PPerv_{G, \zeta}$. The functors $F'_P$ and $F'_{B(M)}$ are compatible with these gradings.
 
  Write $\Vect^{\epsilon}$ for the tensor category of $\ZZ/2\ZZ$-graded vector spaces. Let $\PPerv^{\natural}_{M,G,\zeta}$ be the category of even objects in $\PPerv'_{M, G,\zeta}\otimes \Vect^{\epsilon}$. Let $\PPerv_{G, \zeta}^{\natural}$ be the category of even objects in $\PPerv_{G,\zeta}\otimes \Vect^{\epsilon}$. We get a canonical equivalence of tensor categories $sh: \PPerv^{\natural}_{T,G,\zeta}\,\iso\, \PPerv_{T,G,\zeta}$. The functors $F'_{B(M)}, F'_P, F'_B$ yields tensor functors
\begin{equation}
\label{diag_inclusions_dual_groups}
\PPerv^{\natural}_{G,\zeta} \toup{F^{\natural}_P} \PPerv_{M,G,\zeta}^{\natural} \toup{F^{\natural}_{B(M)}} \PPerv^{\natural}_{T,G,\zeta}
\end{equation}
whose composition is $F^{\natural}_B$. Write $F^{\natural}: \PPerv_{G,\zeta}^{\natural}\to\Vect$ for the functor
$
F_T\comp sh\comp F^{\natural}_B
$. 
By Lemma~\ref{Lm_instead_of_Lm4.2}, $F^{\natural}$ does not annihilate a non-zero object, so it is faithful. By Remark~\ref{Rem_rigid}, $\PPerv^{\natural}_{G,\zeta}$ is a rigid abelian tensor category. Since $F^{\natural}$ is exact and faithful, it is a fibre functor. By (\cite{DM}, Theorem~2.11), $\Aut^{\otimes}(F^{\natural})$ is represented by an affine group scheme $\check{G}_{\zeta}$ over $\Qlb$. We get an equivalence of tensor categories
\begin{equation}
\label{equiv_main_for_checkG} 
\PPerv^{\natural}_{G,\zeta}\,\iso\, \Rep(\check{G}_{\zeta})
\end{equation}
 
 An analog of Remark~\ref{Rem_rigid} holds also for $M$, so $F_T\comp sh\comp F^{\sharp}_{B(M)}: \PPerv^{\natural}_{M,G,\zeta}\to\Vect$ is a fibre functor that yields an affine group scheme $\check{M}_{\zeta}$ and an equivalence of tensor categories $\PPerv^{\natural}_{M,G,\zeta}\,\iso\, \Rep(\check{M}_{\zeta})$. The diagram (\ref{diag_inclusions_dual_groups}) yields homomorphisms $\check{T}_{\zeta}\to \check{M}_{\zeta}\to \check{G}_{\zeta}$. 
 
\subsection{Structure of $\check{G}_{\zeta}$}

\subsubsection{} For $\lambda\in\Lambda^+$ write $\ov{\Gr}^{\lambda}$ for the closure of $\Gr^{\lambda}$ in $\Gr_G$. Let $\ov{\Gra}^{\lambda}_G$ denote the preimage of $\ov{\Gr}^{\lambda}$ in $\Gra_G$.

\begin{Lm} 
\label{Lm_mult_one}
If $\lambda,\mu\in \Lambda^{\sharp,+}_M$ then $\cA^{\lambda+\mu}_{M,\cE}$ appears in $\cA^{\lambda}_{M,\cE}\ast \cA^{\mu}_{M,\cE}$ with multiplicity one.
\end{Lm}
\begin{Prf}
We will give a proof only for $M=G$, the generalization to any $M$ being straightforward. Write $\EE^{\lambda}$ (resp., $\ov{\EE}^{\lambda}$) for the preimage of $\Gra^{\lambda}_G$ (resp., of $\ov{\Gra}^{\lambda}_G$) in $\EE$. As in Section~\ref{Section_Convolution}, we get the convolution map $m^{\lambda,\mu}: \ov{\EE}^{\lambda}\times_{G(\cO)\times \Gm} \ov{\Gra}^{\mu}_G\to \ov{\Gra}^{\lambda+\mu}_G$. Let $W$ be the preimage of $\Gra^{\lambda+\mu}_G$ under $m^{\lambda,\mu}$. Then $m^{\lambda,\mu}$ restricts to an isomorphism $W\to \Gra^{\lambda+\mu}_G$ is an isomorphism, and $W\subset \EE^{\lambda}\times_{G(\cO)\times \Gm} \Gra^{\mu}_G$ is open.
\end{Prf} 
 
\medskip\smallskip

 Write $\Lambda_0=\{\lambda\in\Lambda\mid \<\lambda,\check{\alpha}\>=0\;\mbox{for all}\; \check{\alpha}\in\check{R}\}$. The biggest subgroup in $\Lambda^{\sharp, +}$ is $\Lambda^{\sharp}\cap \Lambda_0$. If $\lambda_1,\ldots, \lambda_r$ generate $\Lambda^{\sharp, +}_M$ as a semi-group then $\oplus_{i=1}^r \cA^{\lambda_i}_{M, \cE}$ is a tensor generator of $\PPerv^{\natural}_{M, G,\zeta}$ in the sense of (\cite{DM}, Proposition~2.20), so $\check{M}_{\zeta}$ is of finite type over $\Qlb$. 
 
 If $\check{M}_{\zeta}$ acts nontrivially on $Z\in\PPerv_{M, G,\zeta}$ then consider the strictly full subcategory of $\PPerv^{\natural}_{M, G,\zeta}$ whose objects are subquotients of $Z^{\oplus m}$, $m\ge 0$. By Lemma~\ref{Lm_mult_one}, this subcategory is not stable under the convolution, so $\check{M}_{\zeta}$ is connected by (\cite{DM}, Corollary~2.22). Since $\PPerv_{M,G,\zeta}$ is semisimple, $\check{M}_{\zeta}$ is reductive by (\cite{DM}, Proposition~2.23). 
 
 By Lemma~\ref{Lm_instead_of_Lm4.2}, for $\lambda\in\Lambda^{\sharp,+}_M$, $w\in W_M$ the weight $w(\lambda)$ of $\check{T}_{\zeta}$ appears in $F^{\natural}(\cA^{\lambda}_{M,\cE})$. So, $\check{T}_{\zeta}$ is closed in $\check{M}_{\zeta}$ by (\cite{DM}, Proposition~2.21). 
 
 For $\nu\in \Lambda^{\sharp,+}_M$ write $V^{\nu}_M$ for the irreducible representation of $\check{M}_{\zeta}$ corresponding to $\cA^{\nu}_{M,\cE}$ via the above equivalence $\PPerv^{\natural}_{M,G,\zeta}\,\iso\, \Rep(\check{M}_{\zeta})$. 
 
\begin{Lm} The torus $\check{T}_{\zeta}$ is maximal in $\check{M}_{\zeta}$. There is a unique Borel subgroup $\check{T}_{\zeta}\subset \check{B}(M)_{\zeta}\subset \check{M}_{\zeta}$ whose set of dominant weights is $\Lambda^{\sharp, +}_M$.
\end{Lm}
\begin{Prf}  First, let us show that for
$\nu_1,\nu_2\in \Lambda^{\sharp,+}_M$ the $\check{T}_{\zeta}$-weight $\nu_1+\nu_2$ appears with multiplicity one in $V^{\nu_1}_M\otimes V^{\nu_2}_M$. 
For $\lambda_1,\lambda_2\in\Lambda$ write $\lambda_1\le\lambda_2$ if $\lambda_2-\lambda_1$ is a sum of some positive coroots for $(G,B)$. By (\cite{MV}, Theorem~3.2) combined with Lemma~\ref{Lm_instead_of_Lm4.2}, if $\nu\in\Lambda^{\sharp}$ appears in $V^{\lambda}_M$ then $w\nu\le \lambda$ for any $w\in W$. By Lemma~\ref{Lm_instead_of_Lm4.2}, the $\check{T}_{\zeta}$-weight $\nu$ appears in $V^{\nu}_M$ with multiplicity one. Our claim follows.

 Let $T'\subset \check{M}_{\zeta}$ be a maximal torus containing $\check{T}_{\zeta}$. By Lemma~\ref{Lm_instead_of_Lm4.2}, for each $\nu\in\Lambda^{\sharp, +}_M$ there is a unique character $\nu'$ of $T'$ such that the two conditions are verified: the composition $\check{T}_{\zeta}\to T'\toup{\nu'}\Gm$ equals $\nu$; the $T'$-weight $\nu'$ appears in $V^{\nu}_M$. The map $\nu\mapsto \nu'$ is a homomoprhism of semigroups, so we can apply (\cite{FL}, Lemma~4.4). This gives a unique Borel subgroup $\check{T}_{\zeta}\subset \check{B}(M)_{\zeta}\subset \check{M}_{\zeta}$ whose set of dominant weights is in bijection with $\Lambda^{\sharp, +}_M$. Since $\nu\mapsto \nu'$ is a bijection between $\Lambda^{\sharp,+}_M$ and the dominant weights of $\check{B}(M)_{\zeta}$, the torus $\check{T}_{\zeta}$ is maximal in $\check{M}_{\zeta}$. 
\end{Prf}
 
\medskip

 For $M=G$ write $\check{B}_{\zeta}=\check{B}(G)_{\zeta}$. So, $\Lambda^{\sharp,+}$ are dominant weights for $(\check{G}_{\zeta}, \check{B}_{\zeta})$. 
If $\lambda\in \Lambda^{\sharp, +}$ lies in the $W$-orbit of $\nu\in \Lambda^{\sharp, +}_M$ then as in Lemma~\ref{Lm_instead_of_Lm4.2} one shows that $\cA^{\nu}_{M,\cE}$ appears in $F^{\natural}_P(\cA^{\lambda}_{\cE})$. By (\cite{DM}, Proposition~2.21) this implies that $\check{M}_{\zeta}$ is closed in $\check{G}_{\zeta}$. 

\subsubsection{Rank one}

Let $M$ be the standard subminimal Levi subgroup of $G$ corresponding to the simple root $\check{\alpha}_i$. Let $j\in J$ be such that $i\in \cI_j$. Let $\check{\Lambda}^{\sharp}=\Hom(\Lambda^{\sharp},\ZZ)$ denote the coweights lattice of $\check{T}_{\zeta}$. Note that $\check{\alpha}_i\in \check{\Lambda}^{\sharp}$. 
Then
$$
\{\check{\nu}\in \check{\Lambda}^{\sharp}\mid \<\lambda,\check{\nu}\>\ge 0\;\mbox{for all}\; \lambda\in \Lambda^{\sharp,+}_M\}
$$
is a $\ZZ_+$-span of a multiple of $\check{\alpha}_i$. So, the group $\check{M}_{\zeta}$ is of semisimple rank 1, and its unique simple coroot is of the form $\check{\alpha}_i/m_i$ for some $m_i\in \QQ$, $m_i>0$. 
 
 Take any $\lambda\in \Lambda^{\sharp, +}_M$ with $\<\lambda, \check{\alpha}_i\>>0$. Write $s_i\in W$ for the simple reflection corresponding to $\check{\alpha}_i$. By Lemma~\ref{Lm_instead_of_Lm4.2}, $F^{\lambda}_M(\cA^{\lambda}_{M,\cE})$ and $F^{s_i(\lambda)}_M(\cA^{\lambda}_{M,\cE})$ do not vanish, so $\lambda-s_i(\lambda)$ is a multiple of the positive root of $\check{M}_{\zeta}$. So, the unique simple root of $\check{M}_{\zeta}$ is $m_i\alpha_i$. It follows that the simple reflection for $(\check{T}_{\zeta}, \check{M}_{\zeta})$ acts on $\Lambda^{\sharp}$ as $\lambda\mapsto \lambda-\<\lambda, \frac{\check{\alpha}_i}{m_i}\> (m_i\alpha_i)=s_i(\lambda)$. We must show that $m_i=\delta_i$. 
 
 By (\cite{MV}, Theorem~3.2) the scheme $\Gr^{\nu}_{B(M)}\cap \Gr^{\lambda}_M$ is non empty if and only if
$$
\nu=\lambda, \lambda-\alpha_i,\lambda-2\alpha_i, \ldots, \lambda-\<\lambda,\check{\alpha}_i\>\alpha_i\, .
$$ 
For $0<k< \<\lambda,\check{\alpha}_i\>$ and $\nu=\lambda-k\alpha_i$ one has
$$
\Gr^{\nu}_{B(M)}\cap \Gr^{\lambda}_M\,\iso\, \Gm\times \AA^{\<\lambda,\check{\alpha}_i\>-k-1}\, .
$$

 Let $M_0$ for the simply-connected cover of the derived group of $M$, Let $T_0$ be the preimage of $T\cap [M,M]$ in $M_0$. Let $\Gr_{M_0}$ denote the affine grassmanian for $M_0$. Let $\Lambda_{M_0}=\ZZ\alpha_i$ denote the coweights lattice of $T_0$. Write $\cL_{M_0}$ for the ample generator of the Picard group of $\Gr_{M_0}$. This is the line bundle with fibre $\det(V_0(\cO): V_0(\cO)^g)$ at $gM_0(\cO)$, where $V_0$ is the standard representation of $M_0$. Let $f_0: \Gr_{M_0}\to \Gr_G$ be the natural map 
 
 For $j'\in J$ the line bundle $f_0^*\cL_{j'}$ is trivial unless $j'=j$, and 
$$
f_0^*\cL_j\,\iso\, \cL_{M_0}^{\frac{\kappa_j(\alpha_i,\alpha_i)}{2}}
$$ 
Besides, the restriction of the line bundle $E_{\beta}/G_{ad}(\cO)$ under $\Gr_{M_0}\toup{f_0}\Gr_G\to \Gr_{G_{ab}}$ is trivial. 

 Assume that $\lambda=a\alpha_i$ with $a>0, a\in\ZZ$ such that $\lambda\in \Lambda^{\sharp}$. Let $\nu=b\alpha_i$ with $b\in\ZZ$ such that $-\lambda<\nu<\lambda$. 

 Write $U\subset M(F)$ for the one-parameter unipotent subgroup corresponding to the affine root space $t^{-a+b}\gg_{\check{\alpha}_i}$. Let $Y$ be the closure of the $U$-orbit through $t^{\nu}M(\cO)$ in $\Gr_{M}$. It is a $T$-stable subscheme $Y\,\iso\, \PP^1$. The $T$-fixed points in $Y$ are $t^{\nu}M(\cO)$ and $t^{-\lambda}M(\cO)$. The natural map $\Gr_{M_0}\to \Gr_M$ induces an isomorphism $\Gr_{M_0}\,\iso\, (\Gr_M^0)_{red}$ at the level of reduced ind-schemes. So, we may consider the restriction of $\cL_{M_0}$ to $Y$, which identifies with $\cO_{\PP^1}(a+b)$. 
  
 The restriction of $\cL_j$ to $\Gr^{\nu}_{B(M)}$ is the constant line bundle with fibre $\Omega_{\bar c}^{\frac{\kappa_{j}(\nu,\nu)}{2}}$. Let $a\in \Omega_{\bar c}^{\frac{\kappa_{j}(\nu,\nu)}{2}}$ be a nonzero element. Viewing it as a section of 
$$
\cL^{\frac{\kappa_{j}(\alpha_i,\alpha_i)}{2}}_{M_0}
$$ 
over $Y$, it will vanish only at $t^{-\lambda}M(\cO)$ with multiplicity $(a+b)\kappa_j(\alpha_i,\alpha_i)/2$. It follows that the shifted  local system $\gt^*_{\nu, B(M)}\cA^{\lambda}_{M,\cE}$ will have the $\Gm$-monodromy 
$$
\zeta^{(a+b)c_j\kappa_{j}(\alpha_i,\alpha_i)/2}
$$ 
This local system is trivial if and only if $(a+b)\bar\kappa(\alpha_i,\alpha_i)/2\in N\ZZ$. We may assume $a\bar\kappa(\alpha_i,\alpha_i)\in 2N\ZZ$. Then the above condition is equivalent to $b\bar\kappa(\alpha_i,\alpha_i)\in 2N\ZZ$. The smallest positive integer $b$ satisfying this condition is $\delta_i$. So, $m_i=\delta_i$. 

\subsubsection{} Let now $M$ be a standard Levi corresponding to a subset $\cI_M\subset \cI$. The semigroup 
$$
\{\check{\nu}\in \check{\Lambda}^{\sharp}\mid \<\lambda,\check{\nu}\>\ge 0\;\mbox{for all}\; \lambda\in \Lambda^{\sharp,+}_M\}
$$
is the $\QQ_+$-closure in $\check{\Lambda}^{\sharp}$ of the $\ZZ_+$-span of positive coroots of $\check{M}_{\zeta}$ with respect to the Borel $\check{B}(M)_{\zeta}$. Since the edges of this convex cone are directed by $\check{\alpha}_i$, $i\in \cI_M$, the simple coroots of $\check{M}_{\zeta}$ are positive rational multiples of $\check{\alpha}_i$, $i\in \cI_M$. Since we know already that $\check{\alpha}_i/\delta_i$, $i\in\cI_M$ are coroots of $\check{M}$, we conclude that the simple coroots of $\check{M}_{\zeta}$ are $\check{\alpha}_i/\delta_i$, $i\in\cI_M$. In turn, this implies that $\check{M}_{\zeta}$ is a Levi subgroup of $\check{G}_{\zeta}$. Finally, we conclude that the Weyl groups of $G$ and of $\check{G}_{\zeta}$ viewed as subgroups of $\Aut(\Lambda^{\sharp})$ are the same.
Theorem~\ref{Th_one} is proved. 

\bigskip\noindent
{\bf Acknowledgements.} We are grateful to V. Lafforgue and M. Finkelberg for fruitfull discussions. The author was supported by the ANR project ANR-13-BS01-0001-01.

\end{document}